\documentclass[12pt]{article}
\usepackage{amsmath}
\usepackage{amssymb}
\usepackage{amsfonts}
\usepackage{amsthm}
\usepackage{enumerate}
\usepackage{mathrsfs}

\def\init{\setcounter{equation}{0}}

\newtheorem{theorem}{Theorem}[section]
\newtheorem{proposition}[theorem]{Proposition}

\newtheorem{definition}[theorem]{Definition}
\newtheorem{corollary}[theorem]{Corollary}

\def\Li{{\operatorname{Li}}}

\def\m{{\operatorname{m}}}

\def\d{{\operatorname{d}}}

\def\F{{\operatorname{F}}}

\def\sgn{{\operatorname{sgn}}}

\def\Re{\hbox{Re}\,}
\def\Im{\hbox{Im}\,}

\numberwithin{equation}{section}
\newenvironment{acknowledgements}{\noindent{\bf Acknowledgements}\bigskip}{}

\begin{document}

\title{Hypergeometric formulas for lattice sums and Mahler measures}
\author{Mathew Rogers\footnote{The author is supported by NSF award DMS-0803107} \\
        \small{\textit{Department of Mathematics, University of
        Illinois}}\\
        \small{\textit{Urbana, IL 61801, USA}}}

\maketitle

\abstract{In this paper, we will prove several formulas relating
generalized hypergeometric functions to lattice sums with four
indices of summation.  These results are related to Boyd's
conjectured identities between Mahler measures and special values
of $L$-series of elliptic curves.}

\section{Introduction}
\label{proofs} \init
In this paper we will prove formulas relating lattice
sums to 
hypergeometric functions. This research was inspired by the work of
Boyd, who used numerical methods to conjecture hundreds of relations
between the $L$-series of elliptic curves and special values of
Mahler's measure \cite{Bo1}. The first example of such an identity
was due to Deninger \cite{De}, who hypothesized that
\begin{equation}\label{deninger's formula}
\m\left(1+y+y^{-1}+z+z^{-1}\right)\stackrel{?}{=}\frac{15}{4\pi^2}L(E,2),
\end{equation}
where $E$ is a conductor $15$ elliptic curve.  As usual,
``$\stackrel{?}{=}$'' denotes a conjectured equality which holds to
at least $50$ decimal places. The Mahler measure of a
$n$-dimensional polynomial is defined in equation \eqref{Mahler
definition}. Boyd observed that since every elliptic curve is
modular, this identity can be translated into a completely explicit
formula:
\begin{equation}\label{deninger's formula after Boyds remark}
\int_{0}^{1}\int_{0}^{1}\log\left|1+2\cos(2\pi t)+2\cos(2\pi
s)\right|\d s\d
t\stackrel{?}{=}\frac{15}{4\pi^2}\sum_{n=1}^{\infty}\frac{a_n}{n^2},
\end{equation}
where
\begin{equation*}
\sum_{n=1}^{\infty}a_n
q^n=q\prod_{n=1}^{\infty}\left(1-q^n\right)\left(1-q^{3n}\right)\left(1-q^{5n}\right)\left(1-q^{15n}\right).
\end{equation*}
If this line of thought is pushed somewhat further (see Theorem
\ref{F(b,c) related to L function}), the following conjecture
arises:
\begin{equation}\label{Deninger intro conjecture latticized}
\begin{split}
\sum_{n=0}^{\infty}&{2n\choose
n}^2\frac{(1/16)^{2n+1}}{2n+1}\\
&\stackrel{?}{=}\frac{540}{\pi^2}\sum_{n_i=-\infty}^{\infty}\frac{(-1)^{n_1+n_2+n_3+n_4}}{\left((6n_1+1)^2+3
(6n_2+1)^2+5(6n_3+1)^2+15(6n_4+1)^2\right)^2}.
\end{split}
\end{equation}
Equation \eqref{Deninger intro conjecture latticized} is an
interesting conjecture, because it relates a complicated lattice sum
to the $_3F_2$ hypergeometric function. Lattice sums have been
extensively studied in physics, since they often arise when
calculating electrostatic potentials of crystal lattices (for
instance see \cite{Gl1}, \cite{Zu} and \cite{FG}).  It is often
difficult to calculate lattice sums numerically, and it is quite
unusual to be able to reduce them to known special functions
\cite{Bor}. The problem of finding a closed form for Madelung's
constant is probably the most famous open problem in this area
\cite{Cr}.
\begin{definition} Let us define $F(a,b,c,d)$ by
\begin{equation*}
\begin{split}
F(a,b,c,d):=&(a+b+c+d)^2\\
&\times\sum_{n_i=-\infty}^{\infty}\frac{(-1)^{n_1+n_2+n_3+n_4}}{\left(
a(6n_1+1)^2+b (6n_2+1)^2+c (6n_3+1)^2+d (6n_4+1)^2\right)^2}.
\end{split}
\end{equation*}
We will also fix the following shorthand notation:
\begin{equation}
F(b,c):=F(1,b,c,b c).
\end{equation}
Finally, the default method of summation will be
$\displaystyle\sum_{n_i=-\infty}^{\infty}=\lim_{v\rightarrow\infty}\sum_{n_1=-v}^{v}\dots\sum_{n_4=-v}^{v}.$
\end{definition}

In this paper, we will prove many new formulas for special values of
$F(a,b,c,d)$.  While we believe that it should be possible to find a
general formula for the function, we have not managed to accomplish
that goal yet.  The second and third sections of the paper summarize
$18$ formulas relating $F(b,c)$ to rational hypergeometric
functions. Virtually all of those identities were extracted directly
from Boyd's tables \cite{Bo1}. For example, Boyd conjectured a
formula for a conductor $20$ elliptic curve, which is equivalent to:
\begin{equation}\label{intro F(1,5) example}
\frac{25}{6\pi^2}F(1,5)\stackrel{?}{=}\sqrt[3]{2} A {_3F_2}\left(\substack{\frac{1}{3},\frac{1}{3},\frac{1}{3}\\
\frac{2}{3},\frac{4}{3}}; \frac{2}{27}\right)
+\sqrt[3]{4} B {_3F_2}\left(\substack{\frac{2}{3},\frac{2}{3},\frac{2}{3}\\
\frac{4}{3},\frac{5}{3}}; \frac{2}{27}\right),
\end{equation}
where $A$ and $B$ are given in terms of gamma functions (see Theorem
\ref{theorem: hypergeometric reduction}).  When more general lattice
sums are considered, hypergeometric functions with irrational
arguments frequently appear. For instance, if
$\phi=\frac{1+\sqrt{5}}{2}$, we have
\begin{equation}\label{intro lattic example}
\frac{225}{32\sqrt{5}\pi^2}F\left(1,5,5,5\right)=\frac{A}{\sqrt[3]{\phi}}{_3F_2}\left(\substack{\frac{1}{3},\frac{1}{3},\frac{1}{3}\\
\frac{2}{3},\frac{4}{3}};\frac{1}{\phi}\right)
+\frac{3B}{\sqrt[3]{\phi^2}}{_3F_2}\left(\substack{\frac{2}{3},\frac{2}{3},\frac{2}{3}\\
\frac{4}{3},\frac{5}{3}};\frac{1}{\phi}\right).
\end{equation}
Formula \eqref{intro lattic example} closely resembles identities
that Forrester and Glasser established for three-dimensional sums
associated with NaCl lattices \cite{FG}.  While it seems likely that
equations \eqref{intro lattic example} and \eqref{intro F(1,5)
example} both arise as special cases of formulas for $F(a,b,c,d)$,
the exact nature of those formulas remains unclear. From our
calculations, we have determined that \eqref{intro lattic example}
is not equivalent to any of Boyd's identities, and it also does not
appear to reduce to a Mahler measure formula.

Based on the computations in this paper, it is probably safe to
conjecture that \textit{many} values of $F(b,c)$ reduce to
generalized hypergeometric functions and Meijer $G$-functions. For
instance, by \eqref{F(1,4) final formula}, we have
\begin{equation*}
\frac{144}{25\pi^2}F(1,4)=m\left(\frac{4}{\theta}\right)+\frac{1}{4\pi^2}\Im\left(G_{3,3}^{3,2}\left(\theta^2\big|\substack{\frac{1}{2},\frac{1}{2},1\\0,0,0}\right)\right),
\end{equation*}
where $m(k)$ is defined in \eqref{m(k) 4F3}, and $\theta\approx
1.93+1.79 i$ is algebraic.  We will briefly describe why
$G$-functions appear in
this context.  
It is often possible to show that lattice sums satisfy
Picard-Fuchs equations with respect to modular parameters.  If the
standard solution to the equation (typically a $_4F_3$ function),
is analytic over the same domain as the lattice sum, then they can
be equated (for an example see \eqref{F(1,2)(x) as a Mahler
measure}). On the other hand, if the two functions have different
domains of analyticity, then the lattice sum will equate to a
piecewise-defined function, which incorporates a second solution
of the differential equation (sometimes a Meijer G-function). One
of our main results is such a formula for $F(1,1,1,x)$, which
holds whenever $x>0$ (see \eqref{F(1,1) hypergeometric form}). We
will also present explicit formulas for $F(1,4)$, $F(2,2)$,
$F(1,1,2,4)$, and $F(1,2,4,4)$. Additionally, we will recover
formulas for $F(1,1)$, $F(1,2)$, and $F(1,3)$, which were first
proved in \cite{RV}.

Therefore, we will briefly outline the approach contained in
Sections \ref{section:hypergeometric reductions}, \ref{section:
F(1,1,2,4) and F(1,2,4,4)}, and \ref{section: F(1,3) section}.
Identities such as those presented in Sections \ref{section : Boyds
conjectures} and \ref{section : summary of hypergeometric formulas}
require the reduction of four-dimensional sums to one-dimensional
sums.  If possible, the first step is a reduction to a
two-dimensional sum. For instance, the left-hand side of
\eqref{intro lattic example} becomes
\begin{equation}\label{intro conductor 32 exam cont}
\begin{split}
F(1,5,5,5)
=16^2\sum_{\substack{n=-\infty\\
k=0}}^{\infty}\frac{(-1)^{n+k}(2k+1)}{\left((6n+1)^2+15(2k+1)^2\right)^2}.
\end{split}
\end{equation}
These sorts of transformations follow from well-known $q$-series
results, and are quite rare. Equation \eqref{intro conductor 32 exam
cont} is a consequence of the following corollary to the Jacobi
triple product:
\begin{equation}\label{intro cm q expansion}
q^2\prod_{n=1}^{\infty}\left(1-q^{3n}\right)\left(1-q^{15n}\right)^3=\sum_{\substack{n=-\infty\\
k=0}}^{\infty}(-1)^{n+k}(2k+1)q^{\frac{15(2k+1)^2+(6n+1)^2}{8}}.
\end{equation}
Notice that equation \eqref{intro cm q expansion} gives an example
of a \textit{lacunary} modular form \cite{Fi2}. The two-dimensional
lattice sums are then evaluated using Ramanujan's theories of
elliptic functions and modular equations. Several of these
calculations are quite involved.

While the main goal of this research was to find formulas for
lattice sums, we have also proved several new relations between
modular forms and Mahler measures. Section \ref{section: more
mahlers} summarizes identities between Mahler measures and Mellin
transforms of non-multiplicative modular forms. Additionally, in
Section \ref{section: conclusion}, we have reformulated a conjecture
concerning a higher Mahler measure.

\section{Summary of Boyd's conjectures for $F(b,c)$}\label{section : Boyds conjectures}\init In this section we
will summarize a variety of explicit formulas relating
four-dimensional lattice sums to Mahler measures of polynomials.
Most of these results are only conjectures, although numerical
calculations can be used to verify them to any degree of accuracy.
Our first step will be to invoke the modularity theorem to find
explicit formulas for $L$-functions of elliptic curves with
conductors $N\in\{11, 14, 15, 20, 24, 27, 32, 36\}$.  The following
theorem is an easy consequence of a paper due to Martin and Ono
\cite{Ono}:

\begin{theorem} Suppose that $E_N$ is an elliptic curve of conductor
$N$, then
\begin{equation}\label{F(b,c) related to L function}
L(E_N,2)=F(b,c)
\end{equation}
for the following values of $N$ and $(b,c)$:
\begin{equation*}
    \begin{tabular}{|c|c|p{6 in}|}
        \hline
        $N$ & $(b,c)$\\
        \hline
        $11$ & $(1,11)$\\
        $14$ & $(2,7)$\\
        $15$ & $(3,5)$\\
        $20$ & $(1,5)$\\
        $24$ & $(2,3)$\\
        $27$ & $(1,3)$\\
        $32$ & $(1,2)$\\
        $36$ & $(1,1)$\\
        \hline
    \end{tabular}
\end{equation*}
\end{theorem}
\begin{proof} We are interested in cases where cusp forms of elliptic
curves equal the product of four eta functions.  Such equalities are
consequences of the modularity theorem. An exhaustive list of all
such cusp forms is provided in \cite{Ono}. By inspection of that
list, the eta product associated with $E_N$ will have the form
\begin{equation*}
g(q):=q\prod_{n=1}^{\infty}\left(1-q^{A n}\right)\left(1-q^{A b
n}\right)\left(1-q^{A c n}\right)\left(1-q^{A b c n}\right),
\end{equation*}
where $(1+b)(1+c)A=24$.  Recalling Euler's pentagonal number theorem
\begin{equation*}
\prod_{n=1}^{\infty}\left(1-q^n\right)=\sum_{n=-\infty}^{\infty}(-1)^n
q^{n (3n+1)/2},
\end{equation*}
this becomes
\begin{equation*}
g(q)=\sum_{n_i=-\infty}^{\infty}(-1)^{n_1+n_2+n_3+n_4}q^{\frac{A
(6n_1+1)^2+A b (6n_2+1)^2+A c (6n_3+1)^2+A b c (6n_4+1)^2}{24}},
\end{equation*}
and it follows immediately that
\begin{equation*}
L(E_N,2)=\frac{24^2}{A^2}\sum_{n_i=-\infty}^{\infty}\frac{(-1)^{n_1+n_2+n_3+n_4}}{\left(
(6n_1+1)^2+b (6n_2+1)^2+c (6n_3+1)^2+b c (6n_4+1)^2\right)^2}.
\end{equation*}
Since $(1+b)(1+c)=24/A$, the theorem follows. $\blacksquare$
\end{proof}

Since we now have expressed several different $L$-values in terms of
$F(b,c)$, it seems logical to list all of the known Mahler measures
which reduce to values of that function.
\begin{definition}  The Mahler measure of
an $n$-dimensional Laurent polynomial, $P(z_1,\dots, z_n)$, is
defined by
\begin{equation}\label{Mahler definition}
\m(P):=\int_{0}^{1}\dots\int_{0}^{1}\log\left|P\left(e^{2\pi i
t_1},\dots, e^{2\pi i t_n}\right)\right| \d t_1\dots \d t_n.
\end{equation}
Furthermore, we will use the following notation:
\begin{align}
m(k):=&\m\left(k+y+y^{-1}+z+z^{-1}\right),\label{m(k) def}\\
n(k):=&\m\left(y^3+z^3+1-k y z\right),\label{n(k) def}\\
g(k):=&\m\left((1+y)(1+z)(y+z)-k y z\right),\label{g(k) def}\\
r(k):=&\m\left((1+y)(1+z)(1+y+z)-k y z\right).\label{r(k) def}
\end{align}
\end{definition}
For convenience we have slightly altered the definitions of $n(k)$,
$g(k)$ and $r(k)$ that appeared in \cite{LR}.  All of the following
examples were either extracted from Boyd's paper \cite{Bo1}, or were
deduced by combining Boyd's conjectures with functional equations in
\cite{LR}. While Boyd's minimal Weierstrass models often do not
coincide with the minimal Weierstrass models in \cite{Ono}, the
elliptic curves are isogenous, and the following results are all
numerically true: \allowdisplaybreaks{
\begin{align}
n(3\sqrt[3]{2})=&\frac{27}{2\pi^2}F(1,1)\\
g(2)=&\frac{9}{2\pi^2}F(1,1)\\
g(-4)=&\frac{18}{\pi^2}F(1,1)\\
m(4 i)=&\frac{16}{\pi^2}F(1,2)\label{m(4i) formula}\\
m(2\sqrt{2})=&\frac{8}{\pi^2}F(1,2)\\
n(-6)=&\frac{81}{4\pi^2}F(1,3)\\
n(\sqrt[3]{2})\stackrel{?}{=}&\frac{25}{6\pi^2}F(1,5)\\
n(\sqrt[3]{32})\stackrel{?}{=}&\frac{40}{3\pi^2}F(1,5)\\
g(-2)\stackrel{?}{=}&\frac{15}{\pi^2}F(1,5)\\
g(4)\stackrel{?}{=}&\frac{10}{\pi^2}F(1,5)\\
r(-1)=&\frac{77}{4\pi^2}F(1,11)\\
m(2)\stackrel{?}{=}&\frac{6}{\pi^2}F(2,3)\\
m(8)\stackrel{?}{=}&\frac{24}{\pi^2}F(2,3)\\
m(3\sqrt{2})\stackrel{?}{=}&\frac{15}{\pi^2}F(2,3)\\
m( i\sqrt{2})\stackrel{?}{=}&\frac{9}{\pi^2}F(2,3)\\
n(-1)\stackrel{?}{=}&\frac{7}{\pi^2}F(2,7)\label{n(-1), F(2,7) conjecture}\\
n(5)\stackrel{?}{=}&\frac{49}{2\pi^2}F(2,7)\\
g(1)\stackrel{?}{=}&\frac{7}{2\pi^2}F(2,7)\\
g(7)\stackrel{?}{=}&\frac{21}{\pi^2}F(2,7)\\
g(-8)\stackrel{?}{=}&\frac{35}{\pi^2}F(2,7)\\
m(1)\stackrel{?}{=}&\frac{15}{4\pi^2}F(3,5)\label{m(1), F(3,5) conjecture}\\
m(3 i)\stackrel{?}{=}&\frac{75}{4\pi^2}F(3,5)\\
m(5)\stackrel{?}{=}&\frac{45}{2\pi^2}F(3,5)\\
m(16)\stackrel{?}{=}&\frac{165}{4\pi^2}F(3,5)
\end{align}
All of the results involving $F(1,1)$, $F(1,2)$, and $F(1,3)$ can be
deduced from Rodriguez-Villegas's paper \cite{RV}. In particular,
those Mahler measures can be written in terms of two-dimensional
Eisenstein-Kronecker series, and then the results follow from
Deuring's theorem.

\section{Summary of rational hypergeometric formulas for
$F(b,c)$}\label{section : summary of hypergeometric formulas}\init

In this section we will translate almost all of the
known Mahler measures for $F(b,c)$ into hypergeometric functions. In
Corollary \ref{final reduction cor} we will also prove that similar
expressions exist for both $F(2,2)$ and $F(1,4)$, even though those
sums are apparently unrelated to the theory of elliptic curves.

\begin{theorem}\label{theorem: hypergeometric reduction} We can express $m(k)$, $n(k)$, and $g(k)$ in terms
of generalized hypergeometric functions for most values of $k$:
\begin{align}
m(k)=&\Re\left(\log(k)-\frac{2}{k^2}{_4F_3}\left(\substack{\frac{3}{2},\frac{3}{2},1,1\\
2,2,2}; \frac{16}{k^2}\right)\right)\label{m(k) 4F3},\\
n(k)=&\Re\left(\log(k)-\frac{2}{k^3}{_4F_3}\left(\substack{\frac{4}{3},\frac{5}{3},1,1\\
2,2,2};\frac{27}{k^3}\right)\right),\label{n(k) 4F3}\\
3g(k)=&\Re\left(\log\left(\frac{(4+k)(k-2)^4}{k^2}\right)-\frac{2k^2}{(4+k)^3}{_4F_3}\left(\substack{\frac{4}{3},\frac{5}{3},1,1\\
2,2,2};\frac{27k^2}{(4+k)^3}\right)\right.\label{g(k) 4F3}\\
&\qquad\left.-\frac{8k}{(k-2)^3}{_4F_3}\left(\substack{\frac{4}{3},\frac{5}{3},1,1\\
2,2,2};\frac{27k}{(k-2)^3}\right)\right)\notag.
\end{align}
Equation \eqref{m(k) 4F3} is valid in $\mathbb{C}\setminus\{0\}$,
while \eqref{n(k) 4F3} and \eqref{g(k) 4F3} are true for $|k|$
sufficiently large, \eqref{g(k) 4F3} also holds in
$\mathbb{R}\setminus[-4,2]$.

In certain cases we can reduce these hypergeometric functions
further. Suppose that $k\in\mathbb{R}\setminus\{0\}$, then
\begin{equation}\label{m(k) 4F3 to 3F2}
\Re\left(\log(k)-\frac{2}{k^2}{_4F_3}\left(\substack{\frac{3}{2},\frac{3}{2},1,1\\
2,2,2}; \frac{16}{k^2}\right)\right)=\Re\left(\frac{|k|}{4}{_3F_2}\left(\substack{\frac{1}{2},\frac{1}{2},\frac{1}{2}\\
1,\frac{3}{2}}; \frac{k^2}{16}\right)\right),
\end{equation}
and
\begin{equation}\label{n(k) 4F3 to 3F2}
\begin{split}
Re\left(\log(k)-\frac{2}{k^3}{_4F_3}\left(\substack{\frac{4}{3},\frac{5}{3},1,1\\
2,2,2};\frac{27}{k^3}\right)\right)=&s(k)\Re\left(A k {_3F_2}\left(\substack{\frac{1}{3},\frac{1}{3},\frac{1}{3}\\
\frac{2}{3},\frac{4}{3}}; \frac{k^3}{27}\right)\right.\\
&\qquad\qquad\left.+B k^2 {_3F_2}\left(\substack{\frac{2}{3},\frac{2}{3},\frac{2}{3}\\
\frac{4}{3},\frac{5}{3}}; \frac{k^3}{27}\right)\right),
\end{split}
\end{equation}
where
$A=\frac{\sqrt[3]{2}\Gamma\left(\frac{1}{6}\right)\Gamma\left(\frac{1}{3}\right)\Gamma\left(\frac{1}{2}\right)}{8\sqrt{3}\pi^2}$,
$B=\frac{\Gamma^3\left(\frac{2}{3}\right)}{16\pi^2}$, and
$s(k)=\frac{1+3\sgn(k)}{4}$.
\end{theorem}
\begin{proof} Equations \eqref{m(k) 4F3} and \eqref{n(k) 4F3} are
due to Rodriguez-Villegas \cite{RV}, while \eqref{g(k) 4F3} was
proved in \cite{LR}.  Kurokawa and Ochiai have examined a version of
\eqref{m(k) 4F3 to 3F2} in \cite{KO}, although it can also be proved
by integrating the following identity:
\begin{equation*}
\Re\left({_2F_1}\left(\substack{\frac{1}{2},\frac{1}{2}\\ 1};
u\right)\right)=\Re\left(\frac{1}{\sqrt{u}}{_2F_1}\left(\substack{\frac{1}{2},\frac{1}{2}\\
1}; \frac{1}{u}\right)\right),
\end{equation*}
which holds for $u\in \left(0,\infty\right)$.  A similar argument
can be used to establish \eqref{n(k) 4F3 to 3F2}. For example, when
$u\in \left(0,\infty\right)$ we can integrate the following
transformation:
\begin{equation*}
\Re\left({_2F_1}\left(\substack{\frac{1}{3},\frac{2}{3}\\ 1};
\frac{1}{u}\right)\right)=\Re\left(\frac{9\Gamma^3\left(2/3\right)}{4\pi^2}
u^{2/3}{_2F_1}\left(\substack{\frac{2}{3},\frac{2}{3}\\
\frac{4}{3}}; u \right)
 +\frac{1}{2} u^{1/3}{_2F_1}\left(\substack{\frac{1}{3},\frac{1}{3}\\
 1};1-u\right)\right),
\end{equation*}
with respect to $u$. $\blacksquare$
\end{proof}

Equations \eqref{m(k) 4F3 to 3F2} and \eqref{n(k) 4F3 to 3F2} will
often allow us to obtain convergent series expansions from divergent
hypergeometric formulas.  For example, applying the results of the
last theorem to conjecture \eqref{m(1), F(3,5) conjecture}, we
obtain formula \eqref{Deninger intro conjecture latticized}. It is
hardly coincidental that \eqref{Deninger intro conjecture
latticized} bears a striking resemblance to a famous formula that
Ramanujan discovered for Catalan's constant \cite{Ad}:
\begin{equation*}
L\left(\chi_{-4},2\right)=\pi\sum_{n=0}^{\infty}{2n\choose
n}^2\frac{\left(1/4\right)^{2n+1}}{2n+1}.
\end{equation*}
Ramanujan's formula follows easily from Boyd's evaluation of the
degenerate Mahler measure $m(4)$.

 The following list summarizes the identities that can be obtained
by reducing the Mahler measures in the previous section to
hypergeometric functions.  Whenever possible, we have used equations
\eqref{m(k) 4F3 to 3F2} and \eqref{n(k) 4F3 to 3F2} to obtain
hypergeometric functions with convergent arguments. Since no such
expression is known for $r(-1)$, we have simply retained that Mahler
measure in our list.  Finally, because several Mahler measures such
as $g(2)$ and $n\left(3\sqrt[3]{2}\right)$ are equivalent, this list
contains fewer entries than we might expect.  Once again, define
\begin{align*}
A:=&\frac{\sqrt[3]{2}\Gamma\left(\frac{1}{6}\right)\Gamma\left(\frac{1}{3}\right)\Gamma\left(\frac{1}{2}\right)}{8\sqrt{3}\pi^2},
&B&:=\frac{\Gamma^3\left(\frac{2}{3}\right)}{16\pi^2},
\end{align*}
then the following results are numerically true:
\begin{align}
\frac{9}{2\pi^2}F(1,1)=&\frac{1}{9}\log(54)-\frac{1}{81}{_4F_3}\left(\substack{\frac{4}{3},\frac{5}{3},1,1\\
2,2,2}; \frac{1}{2}\right),\label{hypergeo F(1,1)}\\
\frac{16}{\pi^2}F(1,2)=&2\log(2)+\frac{1}{8}{_4F_3}\left(\substack{\frac{3}{2},\frac{3}{2},1,1\\
2,2,2}; -\frac{1}{4}\right),\\
\frac{8}{\pi^2}F(1,2)=&\frac{1}{\sqrt{2}}{_3F_2}\left(\substack{\frac{1}{2},\frac{1}{2},\frac{1}{2}\\
1,\frac{3}{2}}; \frac{1}{2}\right),\\
\frac{81}{4\pi^2}F(1,3)=&\log(6)+\frac{1}{108}{_4F_3}\left(\substack{\frac{4}{3},\frac{5}{3},1,1\\
2,2,2}; -\frac{1}{8}\right),\label{hypergeo F(1,3)}\\
\frac{25}{6\pi^2}F(1,5)\stackrel{?}{=}&\sqrt[3]{2} A {_3F_2}\left(\substack{\frac{1}{3},\frac{1}{3},\frac{1}{3}\\
\frac{2}{3},\frac{4}{3}}; \frac{2}{27}\right)
+\sqrt[3]{4} B {_3F_2}\left(\substack{\frac{2}{3},\frac{2}{3},\frac{2}{3}\\
\frac{4}{3},\frac{5}{3}}; \frac{2}{27}\right),\\
\frac{40}{3\pi^2}F(1,5)\stackrel{?}{=}&\frac{5}{3}\log(2)-\frac{1}{16}{_4F_3}\left(\substack{\frac{4}{3},\frac{5}{3},1,1\\
2,2,2}; \frac{27}{32}\right),\\
\frac{77}{4\pi^2}F(1,11)=&r(-1),\label{hypergeo F(1,11)}\\
\frac{6}{\pi^2}F(2,3)\stackrel{?}{=}&\frac{1}{2}{_3F_2}\left(\substack{\frac{1}{2},\frac{1}{2},\frac{1}{2}\\
1,\frac{3}{2}}; \frac{1}{4}\right),\label{hypergeo F(2,3) first 0}\\
\frac{24}{\pi^2}F(2,3)\stackrel{?}{=}&3\log(2)-\frac{1}{32}{_4F_3}\left(\substack{\frac{3}{2},\frac{3}{2},1,1\\
2,2,2}; \frac{1}{4}\right),\label{hypergeo F(2,3) first}\\
\frac{15}{\pi^2}F(2,3)\stackrel{?}{=}&\frac{1}{2}\log(18)-\frac{1}{9}{_4F_3}\left(\substack{\frac{3}{2},\frac{3}{2},1,1\\
2,2,2}; \frac{8}{9}\right),\\
\frac{9}{\pi^2}F(2,3)\stackrel{?}{=}&\frac{1}{2}\log(2)+{_4F_3}\left(\substack{\frac{3}{2},\frac{3}{2},1,1\\{2,2,2}},-8\right),\label{hypergeo F(2,3) last}\\
\frac{7}{\pi^2}F(2,7)\stackrel{?}{=}&\frac{A}{2} {_3F_2}\left(\substack{\frac{1}{3},\frac{1}{3},\frac{1}{3}\\
\frac{2}{3},\frac{4}{3}}; -\frac{1}{27}\right)
-\frac{B}{2} {_3F_2}\left(\substack{\frac{2}{3},\frac{2}{3},\frac{2}{3}\\
\frac{4}{3},\frac{5}{3}}; -\frac{1}{27}\right),\\
\frac{49}{2\pi^2}F(2,7)\stackrel{?}{=}&\log(5)-\frac{2}{125}{_4F_3}\left(\substack{\frac{4}{3},\frac{5}{3},1,1\\
2,2,2}; \frac{27}{125}\right),\\
\frac{21}{\pi^2}F(2,7)\stackrel{?}{=}&g(7),\\
\frac{15}{4\pi^2}F(3,5)\stackrel{?}{=}&\frac{1}{4}{_3F_2}\left(\substack{\frac{1}{2},\frac{1}{2},\frac{1}{2}\\
1,\frac{3}{2}}; \frac{1}{16}\right),\label{f(3,5) first}\\
\frac{45}{2\pi^2}F(3,5)\stackrel{?}{=}&\log(5)-\frac{2}{25}{_4F_3}\left(\substack{\frac{3}{2},\frac{3}{2},1,1\\
2,2,2}; \frac{16}{25}\right),\\
\frac{165}{4\pi^2}F(3,5)\stackrel{?}{=}&4\log(2)-\frac{1}{128}{_4F_3}\left(\substack{\frac{3}{2},\frac{3}{2},1,1\\
2,2,2}; \frac{1}{16}\right),\\
\frac{75}{4\pi^2}F(3,5)\stackrel{?}{=}&\log(3)+\frac{2}{9}{_4F_3}\left(\substack{\frac{3}{2},\frac{3}{2},1,1\\
2,2,2}; -\frac{16}{9}\right).\label{hypergeo very last}
\end{align}
While most of these formulas remain unproven, a variety of partial
results exist.  For instance, identities \eqref{hypergeo F(2,3)
first 0}} through \eqref{hypergeo F(2,3) last} are equivalent to one
another \cite{LR}, formulas \eqref{f(3,5) first} through
\eqref{hypergeo very last} are equivalent to one another \cite{Ml2},
and formulas \eqref{hypergeo F(1,1)} through \eqref{hypergeo F(1,3)}
follow from \cite{RV}.  Additionally, Mellit and Brunault have given
$K$-theoretic proofs of the formulas for $F(2,7)$ and $F(1,11)$ (see
\cite{Me} and \cite{Br}).

\section{Proofs of formulas for $F(1,1)$, $F(1,2)$, $F(1,4)$, $F(2,2)$, and $F(1,1,1,x)$}\label{section:hypergeometric reductions}\init

In the previous section we translated many of Boyd's conjectures
into explicit identities between hypergeometric functions and
lattice sums.  This approach has two essential consequences. Not
only does it eliminate any obvious connection with elliptic curves,
but it also allows for the construction of proofs based upon series
manipulation.  In this section we will discuss the cases that occur
when $(b,c)\in\{(1,1),(1,2),(1,4),(2,2)\}$.  We will rely heavily on
the $q$-series theorems contained in Ramanujan's notebooks (see
\cite{Be3} and \cite{Be5}).

\begin{definition} Let us recall Ramanujan's $q$-series notation:
\begin{align*}
\varphi(q):=&\sum_{n=-\infty}^{\infty}q^{n^2},&
\psi(q):=&\sum_{n=0}^{\infty}q^{\frac{n(n+1)}{2}},\\
 f(-q):=&\prod_{n=1}^{\infty}(1-q^n),&
(x;q)_{\infty}:=&\prod_{n=0}^{\infty}\left(1-x q^n\right).
\end{align*}
When convenient, we will also employ the following notation:
\begin{equation*}
e_j:=q^{j/24}\prod_{n=1}^{\infty}\left(1-q^{j n}\right).
\end{equation*}
\end{definition}

Proposition \ref{lemma: cm reduction} reduces the aforementioned
cases of $F(b,c)$ to two-dimensional sums.  Such identities exist
because various eta-quotients can be written in terms theta
functions. Euler's pentagonal number formula is probably the
simplest such identity:
\begin{equation*}
e_j=\sum_{n=-\infty}^{\infty}(-1)^n q^{\frac{j(6n+1)^2}{24}}.
\end{equation*}
Unfortunately, similar formulas are not known (and probably do not
exist) for $e_1^2$, $e_1e_2$, or $e_1e_3$ \cite{Fi}. This fact
represents the main obstruction to proving Boyd's conjectures for
$F(1,5)$, $F(2,3)$, $F(2,7)$, and $F(3,5)$.

\begin{definition} We will use the following notation:
\begin{align}
F_{(1,2)}(x):=&\sum_{\substack{n=-\infty\\ k=0}}^{\infty}\frac{(-1)^{n+k}(2k+1)}{\left((2k+1)^2+x^2(2n)^2\right)^2},\label{F(1,2) 2D sum}\\
F_{(1,4)}(x):=&25\sum_{n,k=-\infty}^{\infty}\frac{(-1)^n
(3k+1)}{\left(4(3k+1)^2+x^2(6n+1)^2\right)^2},\label{F(1,4) 2D sum}\\
F_{(2,2)}(x):=&9\sum_{n,k=0}^{\infty}\frac{(-1)^{\frac{n(n+1)}{2}+k}(2k+1)}{\left(2(2k+1)^2+x^2(2n+1)^2\right)^2}.\label{F(2,2)
2D sum}
\end{align}
\end{definition}

In the next proposition we will show that each of these functions
equals a value of $F(b,c)$ when $x\rightarrow 1$. Furthermore, we
will also demonstrate that a two-dimensional series exists for
$F(1,1,1,x)$.

\begin{proposition}\label{lemma: cm reduction} Suppose that
$(b,c)\in\{(1,2),(1,4),(2,2)\}$, then
\begin{equation}\label{4d to 2d sum reduction}
F_{(b,c)}(1)=F(b,c).
\end{equation}
Furthermore, the following series expansion is true for any $x>0$:
\begin{equation}
\frac{F\left(1,1,1,x\right)}{\left(3+x\right)^2}=\sum_{\substack{n=-\infty\\
k=0}}^{\infty}\frac{(-1)^{n+k}(2k+1)}{\left(3(2k+1)^2+x
(6n+1)^2\right)^2}.\label{F(1,1) 2D sum}
\end{equation}
\end{proposition}
\begin{proof}
First notice that $F(a,b,c,d)$ has the following integral
representation for all positive values of $a$, $b$, $c$, and $d$:
\begin{equation*}
\frac{24^2 F(a,b,c,d)}{(a+b+c+d)^2} =\int_{0}^{1}\int_{0}^{q_1}e_a
e_b e_c e_d\frac{\d q}{q}\frac{\d q_1}{q_1}.
\end{equation*}
Taking note of the following identities:
\begin{align*}
e_x e_1^3=&(e_x)e_1^3,\\
e_1^2 e_2^2=&\left(\frac{e_1^2}{e_2}\right)e_2^3,\\
e_1^2 e_4^2=&\left(\frac{e_1^2e_4^2}{e_2}\right)e_2,\\
e_1 e_2^2 e_4=&\left(\frac{e_1e_4}{e_2}\right)e_2^3,
\end{align*}
and then employing well known series expansions:
\begin{align}
\frac{e_1^2}{e_2}=&\sum_{n=-\infty}^{\infty}(-1)^n q^{n^2},\\
\frac{e_1e_4}{e_2}=&\sum_{n=0}^{\infty}(-1)^{\frac{n(n+1)}{2}}q^{\frac{(2n+1)^2}{8}},\label{weight 1/2 series}\\
e_j=&\sum_{n=-\infty}^{\infty}(-1)^n q^{\frac{j(6n+1)^2}{24}},\\
e_j^3=&\sum_{n=0}^{\infty}(-1)^n (2n+1)q^{\frac{j(2n+1)^2}{8}},\\
\frac{e_1^2e_4^2}{e_2}=&\sum_{n=-\infty}^{\infty}(3n+1)q^{\frac{(3n+1)^2}{3}}\label{weight
3/2 series},
\end{align}
we recover equations \eqref{F(1,1) 2D sum} and \eqref{4d to 2d sum
reduction} in every case. $\blacksquare$
\end{proof}

We will use the next proposition to reduce each of the
two-dimensional sums to a $q$-series.  Then, in Theorem
\ref{hypergeometric theorem}, we will reduce each $q$-series to an
integral of hypergeometric functions. In certain special cases those
integrals translate into identities involving Mahler measures.

\begin{proposition}\label{q-series proposition} Let $\chi_{-3}(k)$ and $\chi_{-4}(k)$ denote Legendre
symbols modulo three and four, and assume that $x>0$.

If $q=e^{-\pi x/\sqrt{12}}$ and $\omega=e^{\pi i/6}$, then
\begin{align}\label{F(1,1)(x) q series sum}
\frac{F\left(1,1,1,x^2\right)}{\left(3+x^2\right)^2}&=\frac{\pi^2}{72x}\sum_{k=1}^{\infty}k\chi_{-4}(k)
\log\left|\frac{1+\omega q^{k}}{1-\omega q^{k}}\right|.
\end{align}
If $q=e^{-\pi x}$, then
\begin{equation}\label{F(1,2)(x) q series sum}
F_{(1,2)}(x)=-\frac{\pi^2}{32x}\left(\log(q)+4\sum_{k=1}^{\infty}k\chi_{-4}(k)\log\left(1+q^{k}\right)\right).
\end{equation}
If $q=e^{-\pi x/3}$ and $\omega=e^{\pi  i/6}$, then
\begin{align}\label{F(1,4)(x) q-series sum}
F_{(1,4)}(x)&=\frac{25\pi^2}{72x}\sum_{k=1}^{\infty}k\chi_{-3}(k)\log\left|\frac{1+\omega
q^k}{1-\omega q^k}\right|.
\end{align}
If $q=e^{-\pi x/\sqrt{8}}$ and $\omega=e^{\pi i/4}$, then
\begin{align}\label{F(2,2)(x) q-series sum}
F_{(2,2)}(x)&=\frac{9\pi^2}{32x}\sum_{k=1}^{\infty}
k\chi_{-4}(k)\log\left|\frac{1+\omega q^{k}}{1-\omega q^{k}}\right|.
\end{align}
\end{proposition}
\begin{proof} All of the proofs are very similar, so we
will only prove \eqref{F(1,1)(x) q series sum} in detail. First
notice that
\begin{align*}
\frac{F(1,1,1,x^2)}{\left(3+x^2\right)^2}=&\frac{1}{144}\sum_{\substack{n=-\infty\\
k=0}}^{\infty}\frac{(-1)^{n+k}(2k+1)}{\left((k+1/2)^2+\frac{x^2}{12}
(6n+1)^2\right)^2},\\
=&\frac{\pi^2}{144}\int_{0}^{\infty}u\left(\sum_{
k=0}^{\infty}(-1)^{k}(2k+1)e^{-\pi(k+1/2)^2
u}\right)\\
&\qquad\qquad\times\left(\sum_{n=-\infty}^{\infty}(-1)^n e^{-\pi
(6n+1)^2 x^2 u/12}\right)\d u.
\end{align*}
Next, by the involution for the weight $3/2$ theta function:
\begin{equation*}
\sum_{k=0}^{\infty}(-1)^k (2k+1)e^{-\pi(k+1/2)^2
u}=\frac{1}{u^{3/2}}\sum_{k=0}^{\infty}(-1)^k (2k+1)e^{-\pi
(k+1/2)^2\frac{1}{u}},
\end{equation*}
this becomes
\begin{align*}
\frac{F(1,1,1,x^2)}{\left(3+x^2\right)^2}=&\frac{\pi^2}{144}\sum_{\substack{n=-\infty\\
k=0}}^{\infty}(-1)^{n+k}(2k+1)\int_{0}^{\infty}u^{-1/2}e^{-\frac{\pi(k+1/2)^2}{
u}- \frac{\pi(6n+1)^2 x^2}{12}u}\d u.
\end{align*}
The following integral holds whenever $A,B\in\mathbb{R}$ and
$|A|\not=0$:
\begin{equation*}
\int_{0}^{\infty}u^{-1/2}e^{-\pi \left(A^2 u+\frac{B^2}{u}\right)}\d
u=\frac{e^{-2\pi|AB|}}{|A|},
\end{equation*}
and therefore
\begin{align*}
\frac{F(1,1,1,x^2)}{\left(3+x^2\right)^2}=&\frac{\pi^2}{24\sqrt{3}x}\sum_{\substack{n=-\infty\\
k=0}}^{\infty}(-1)^{n+k}\frac{(2k+1)}{|6n+1|}e^{-\pi(2k+1)|6n+1|x/\sqrt{12}}\\
=&\frac{\pi^2}{24\sqrt{3}x}\sum_{k=1}^{\infty}k\chi_{-4}(k)\sum_{n=-\infty}^{\infty}\frac{(-1)^{n}}{|6n+1|}q^{k|6n+1|}\\
=&\frac{\pi^2}{72x}\sum_{k=1}^{\infty}k\chi_{-4}(k)\log\left|\frac{1+\omega
q^{k}}{1-\omega q^{k}}\right|,
\end{align*}
where $\omega=e^{\pi i/6}$ and $q=e^{-\pi x/\sqrt{12}}$.
$\blacksquare$
\end{proof}

At this point, a hypergeometric formula for $F_{(1,2)}(x)$ can be
recovered.  By combining equation \eqref{F(1,2)(x) q series sum}
with formulas (2-9) and (2-16) in \cite{LR}, it is possible to show
that if $q=e^{-\pi x}$, then
\begin{equation}\label{F(1,2)(x) as a Mahler measure}
F_{(1,2)}(x)=\frac{\pi^2}{16x}m\left(\frac{ i
f^4(-q)}{\sqrt{q}f^4\left(-q^4\right)}\right).
\end{equation}
In the next section we will use values of class invariants to deduce
explicit examples from \eqref{F(1,2)(x) as a Mahler measure}.
Unfortunately, we will require another theorem to obtain useful
results on the other lattice sums.

\begin{theorem}\label{q-series integrals theorem} In this theorem we will
always assume that $x>0$.  If $q=e^{-\pi x/\sqrt{12}}$, then
\begin{align}\label{F(1,1)(x) q-series integral}
\frac{F(1,1,1,x^2)}{\left(3+x^2\right)^2}&=\frac{\pi^2}{24\sqrt{3}x}\Im\left(\int_{0}^{
i q}\frac{f^9\left(-u^3\right)}{f^3\left(-u\right)}\d u\right).
\end{align}
If $q=e^{-\pi x}$, then
\begin{equation}\label{F(1,2)(x) q-series integral}
F_{(1,2)}(x)=\frac{\pi^3}{32}-\frac{\pi^2}{16x}\int_{0}^{q}\frac{\varphi^2(-u)\varphi^4(u)-1}{u}\d
u.
\end{equation}
If $q=e^{-\pi x/3}$ and $\rho=e^{2\pi i/3}$, then
\begin{align}\label{F(1,4)(x) q-series integral}
F_{(1,4)}(x)=\frac{25\pi^2}{36 x}\Im\left(\int_{0}^{\rho
q}\varphi^2(u)\psi^4\left(u^2\right)\d u\right).
\end{align}
If $q=e^{-\pi x/\sqrt{8}}$, then
\begin{align}\label{F(2,2)(x) q-series integral}
F_{(2,2)}(x)&=\frac{9\pi^2}{32\sqrt{2} x}\int_{0}^{q}
\varphi(-u^2)\varphi(u^4)\left(3\psi^4(-u^2)-\psi^4(u^2)\right)\d u.
\end{align}
\end{theorem}
\begin{proof}  Equations \eqref{F(1,2)(x) q-series integral} and \eqref{F(1,4)(x)
q-series integral} have similar proofs, so we will only prove the
latter identity. Notice that \eqref{F(1,4)(x) q-series sum} can be
rearranged to obtain
\begin{align*}
\frac{72x}{25\pi^2}F_{(1,4)}(x)=&\Re\left(\int_{0}^{q}
\sum_{k=1}^{\infty}k^2\chi_{-3}(k)\frac{2\omega u^k}{1-\omega^2
u^{2k}}\frac{\d u}{u}\right)\\
=&\sqrt{3}\int_{0}^{q}\sum_{k=1}^{\infty}k^2\chi_{-3}(k)\left(\frac{u^k-u^{5k}}{1+u^{6k}}\right)\frac{\d
u}{u}\\
=&\Im\left(2\int_{0}^{\rho
q}\sum_{k=1}^{\infty}k^2\left(\frac{u^k-u^{3k}+u^{5k}}{1+u^{6k}}\right)\frac{\d
u}{u}\right)\\
=&\Im\left(2\int_{0}^{\rho q}\sum_{k=1}^{\infty}\frac{k^2
u^k}{1+u^{2k}}\frac{\d u}{u}\right),
\end{align*}
where $\rho=e^{2\pi i/3}$. Combining entries 10.1, 11.3, and 17.2 in
Chapter 17 of \cite{Be3}, we deduce that for $|u|<1$:
\begin{equation*}
\sum_{k=1}^{\infty}\frac{k^2
u^k}{1+u^{2k}}=u\varphi^2(u)\psi^4\left(u^2\right),
\end{equation*}
which completes the proof of \eqref{F(1,4)(x) q-series integral}.

The proofs of equations \eqref{F(1,1)(x) q-series integral} and
\eqref{F(2,2)(x) q-series integral} will require the following
formula:
\begin{equation}\label{fine identity}
\Im\left(g\left( i
u,t\right)\right)=\frac{1}{t}\sum_{k=0}^{\infty}(-1)^k\frac{(2k+1)^2
u^{2k+1}}{1+u^{2(2k+1)}}\left(t^{2k+1}-t^{-(2k+1)}\right),
\end{equation}
where
\begin{equation*}
g(u,t)=\frac{(u;u)_{\infty}^6\left(t^{-2}u;u\right)_{\infty}
\left(t^2;u\right)_{\infty}}{\left(t^{-1}u;u\right)_\infty^4
\left(t;u\right)_{\infty}^4}.
\end{equation*}
Equation \eqref{fine identity} is a direct consequence of identity
(14.2.9) in \cite{BA}, which follows from product expansions for the
Weierstrass $\wp$-function \cite{Ca}.

Rearranging equation \eqref{F(1,1)(x) q series sum} we have
\begin{align*}
\frac{72x}{\pi^2}\frac{F(1,1,1,x^2)}{\left(3+x^2\right)^2}=&\Re\left(\int_{0}^{q}
\sum_{k=0}^{\infty}(-1)^k (2k+1)^2\frac{2\omega u^{2k+1}}{1-\omega^2
u^{2(2k+1)}}\frac{\d u}{u}\right)\\
=&\sqrt{3}\int_{0}^{q} \sum_{k=0}^{\infty}(-1)^k
(2k+1)^2\frac{u^{2k+1}-u^{5(2k+1)}}{1+u^{6(2k+1)}}\frac{\d u}{u}.
\end{align*}
Applying \eqref{fine identity} after letting $u\rightarrow u^3$ and
$t\rightarrow u^{-2}$, transforms this last integral into
\begin{equation*}
\frac{24\sqrt{3}x}{\pi^2}\frac{F(1,1,1,x^2)}{\left(3+x^2\right)^2}=\Im\left(\int_{0}^{
i
q}\frac{f^5\left(-u^2\right)f^4\left(-u^3\right)f\left(-u^6\right)}{f^4(-u)}\d
u\right).
\end{equation*}
By the following eta function identity:
\begin{equation}\label{eta function identity in signature three}
\frac{f^5\left(-u^2\right)f^4\left(-u^3\right)f\left(-u^6\right)}{f^4(-u)}=\frac{f^9\left(-u^3\right)}{f^3\left(-u\right)}
+u\frac{f^9\left(-u^6\right)}{f^3\left(-u^2\right)},
\end{equation}
this becomes
\begin{align*}
\frac{24\sqrt{3}x}{\pi^2}\frac{F(1,1,1,x^2)}{\left(3+x^2\right)^2}=&\Im\left(\int_{0}^{
i q}\frac{f^9\left(-u^3\right)}{f^3\left(-u\right)}\d
u\right)-\Im\left(\int_{0}^{q}
u\frac{f^9\left(u^6\right)}{f^3\left(u^2\right)}\d u\right)\\
=&\Im\left(\int_{0}^{ i
q}\frac{f^9\left(-u^3\right)}{f^3\left(-u\right)}\d u\right)-0,
\end{align*}
which completes the proof of \eqref{F(1,1)(x) q-series integral}.
Although we will not elaborate on the proof of \eqref{eta function
identity in signature three} here, it suffices to say that it
follows from algebraic transformations for the hypergeometric
function.

The proof of \eqref{F(2,2)(x) q-series integral} follows the same
lines, but requires a few extra steps.  Proceeding as before, we
find that
\begin{align*}
\frac{32x}{9\pi^2}F_{(2,2)}(x)=&\Re\left(\int_{0}^{q}
\sum_{k=0}^{\infty}(-1)^k (2k+1)^2\frac{2\omega u^{2k+1}}{1-\omega^2
u^{2(2k+1)}}\frac{\d u}{u}\right)\\
 =&\int_{0}^{q}
\sum_{k=0}^{\infty}(-1)^k (2k+1)^2\frac{u^{2k+1}-u^{3(2k+1)}}{1+
u^{4(2k+1)}}\frac{\d u}{u}.
\end{align*}
Applying equation \eqref{fine identity} after letting $u\rightarrow
u^2$ and $t\rightarrow u^{-1}$, this becomes
\begin{align*}
\frac{32x}{9\pi^2}F_{(2,2)}(x)=\Re\left(2\int_{0}^{\omega q}
\frac{f^6\left(-u^2\right)f\left(-u^8\right)}{f\left(-u^4\right)}\frac{1}{\left(\omega
u, u^2\right)_{\infty}^4\left(\bar{\omega}u,u^2\right)_{\infty}^4}\d
u\right),
\end{align*}
where $\omega=e^{\pi i/4}$.  For brevity of notation let us define a
new function
\begin{align}
g(u):=&\left(\omega u,
u^2\right)_{\infty}\left(\bar{\omega}u,u^2\right)_{\infty}\notag\\
=&\prod_{n=0}^{\infty}\left(1-\sqrt{2}
u^{2n+1}+u^{2(2n+1)}\right)\label{g(u) product}.
\end{align}
Since \eqref{F(2,2)(x) q-series sum} is odd with respect to $q$, our
integral can be transformed into
\begin{align}\label{F(2,2)(x) qintegral intermed}
\frac{32x}{9\pi^2}F_{(2,2)}(x)=\Re\left(\int_{0}^{\omega q}
\frac{f^6\left(-u^2\right)f\left(-u^8\right)}{f\left(-u^4\right)}\frac{g^4(u)+g^4(-u)}{g^4(u)g^4(-u)}\d
u\right).
\end{align}
Next we will reduce
$\left(g^4(u)+g^4(-u)\right)/\left(g(u)g(-u)\right)^4$ to theta
functions.  Observe by equation \eqref{g(u) product} that
\begin{align}\label{g(u)g(-u) product}
g(u)g(-u)=&\prod_{n=0}^{\infty}\left(1+u^{4(2n+1)}\right)
=\frac{\varphi\left(-u^8\right)}{f\left(-u^4\right)}.
\end{align}
With two applications of the Jacobi triple product \cite{Be3}, we
also have
\begin{align}
g(u)+g(-u)=&\frac{1}{f\left(-u^2\right)}\left(\sum_{n=-\infty}^{\infty}(-\omega)^n
u^{n^2}+\sum_{n=-\infty}^{\infty}\omega^n u^{n^2}\right)\notag\\
=&\frac{2}{f\left(-u^2\right)}\sum_{n=-\infty}^{\infty} (-1)^{n}
u^{16n^2}\notag\\
=&\frac{2\varphi\left(-u^{16}\right)}{f\left(-u^2\right)}
\label{g(u)+g(-u) sum}.
\end{align}
So finally, combining \eqref{g(u)g(-u) product} and
\eqref{g(u)+g(-u) sum}, we find that
\begin{equation*}
\frac{g^4(u)+g^4(-u)}{g^4(u)g^4(-u)}
=2\frac{f^4\left(-u^4\right)}{f^4\left(-u^2\right)}
\left(8\frac{\varphi^4\left(-u^{16}\right)}{\varphi^4\left(-u^8\right)}-8\frac{\varphi^2\left(-u^{16}\right)\varphi\left(-u^2\right)}{\varphi^3\left(-u^8\right)}+\frac{\varphi^2\left(-u^2\right)}{\varphi^2\left(-u^8\right)}\right).
\end{equation*}
Recalling that
$\varphi^2\left(-q^{16}\right)=\varphi\left(-q^{8}\right)\varphi\left(q^8\right)$,
this becomes
\begin{align*}
\frac{g^4(u)+g^4(-u)}{g^4(u)g^4(-u)}
=&2\frac{f^4\left(-u^4\right)}{f^4\left(-u^2\right)}
\left(\frac{8\varphi^2\left(u^{8}\right)-8\varphi\left(u^8\right)\varphi\left(-u^2\right)+\varphi^2\left(-u^2\right)}{\varphi^2\left(-u^8\right)}\right),\\
=&2\frac{f^5\left(-u^4\right)}{f^6\left(-u^2\right)}
\left(\frac{8\varphi^2\left(u^{8}\right)\varphi\left(-u^2\right)-8\varphi\left(u^8\right)\varphi^2\left(-u^2\right)+\varphi^3\left(-u^2\right)}{\varphi^2\left(-u^8\right)}\right),
\end{align*}
and therefore \eqref{F(2,2)(x) qintegral intermed} simplifies to

\begin{equation*}
\begin{split}
\frac{32x}{9\pi^2}F_{(2,2)}(x)=&\Re\left(2\int_{0}^{\omega q}
f^4\left(-u^4\right)f\left(-u^8\right)\right.\\
&\qquad\qquad\left.\times\frac{8\varphi^2\left(u^8\right)\varphi\left(-u^2\right)-8\varphi\left(u^8\right)\varphi^2\left(-u^2\right)+\varphi^3\left(-u^2\right)}{\varphi^2\left(-u^8\right)}\d
u\right).
\end{split}
\end{equation*}
Next, let $u\rightarrow \omega u$ to obtain
\begin{equation*}
\begin{split}
=&2\int_{0}^{q}
\frac{f^4\left(u^4\right)f\left(-u^8\right)}{\varphi^2\left(-u^8\right)}\\
&\qquad\times\Re\left(\omega\left(8\varphi^2\left(u^8\right)\varphi\left(-
i u^2\right)-8\varphi\left(u^8\right)\varphi^2\left(- i
u^2\right)+\varphi^3\left(- i u^2\right)\right)\right)\d u.
\end{split}
\end{equation*}
If we recall that $\varphi\left(- i
u^2\right)=\varphi\left(u^8\right)-2 i u^2\psi\left(u^{16}\right)$,
then we are left with
\begin{align*}
\frac{32x}{9\pi^2}F_{(2,2)}(x) =&\sqrt{2}\int_{0}^{q}
\frac{f^4\left(u^4\right)f\left(-u^8\right)}{\varphi^2\left(-u^8\right)}\\
&\qquad\times\left(\left(\varphi\left(u^8\right)-2u^2\psi\left(u^{16}\right)\right)^3
-4u^2\varphi\left(u^8\right)\psi\left(u^{16}\right)\left(\varphi\left(u^8\right)-2u^2\psi\left(u^{16}\right)\right)\right)\d
u.
\end{align*}
In order to simplify this last formula, we will freely apply theta
function identities on pages 34 and 40 of \cite{Be3}.  Therefore, we
find that
\begin{align*}
=&\sqrt{2}\int_{0}^{q}\frac{f^4\left(u^4\right)f\left(-u^8\right)}{\varphi^2\left(-u^8\right)}\left(\varphi^3\left(-u^2\right)-4u^2\varphi\left(u^8\right)\psi\left(u^{16}\right)\varphi\left(-u^2\right)\right)\d
u\\
=&\sqrt{2}\int_{0}^{q}\frac{f^4\left(u^4\right)f\left(-u^8\right)}{\varphi^2\left(-u^8\right)}\varphi\left(-u^2\right)\left(\varphi^2\left(-u^2\right)-4u^2\psi^2\left(u^{8}\right)\right)\d
u\\
=&\frac{1}{\sqrt{2}}\int_{0}^{q}\frac{f^4\left(u^4\right)f\left(-u^8\right)}{\varphi^2\left(-u^8\right)}\varphi\left(-u^2\right)\left(3\varphi^2\left(-u^2\right)-\varphi^2\left(u^{2}\right)\right)\d
u\\
=&\frac{1}{\sqrt{2}}\int_{0}^{q}\frac{f^4\left(u^4\right)f\left(-u^8\right)}{\varphi^2\left(-u^8\right)\psi^2\left(u^4\right)}
\varphi\left(-u^2\right)\left(3\psi^4\left(-u^2\right)-\psi^4\left(u^{2}\right)\right)\d
u\\
=&\frac{1}{\sqrt{2}}\int_{0}^{q}\varphi\left(u^4\right)
\varphi\left(-u^2\right)\left(3\psi^4\left(-u^2\right)-\psi^4\left(u^{2}\right)\right)\d
u,
\end{align*}
which completes the proof of \eqref{F(2,2)(x) q-series integral}.
$\blacksquare$
\end{proof}

The next theorem requires the signature-three theta functions.
Recall that if $\omega=e^{2\pi i/3}$, then the signature-three theta
functions are defined by:
\begin{align*}
a(q):=&\sum_{n,m=-\infty}^{\infty}q^{m^2+m n+n^2},\\
b(q):=&\sum_{n,m=-\infty}^{\infty}\omega^{m-n}q^{m^2+m n+n^2},\\
c(q):=&\sum_{n,m=-\infty}^{\infty}q^{(m+1/3)^2+(m+1/3)(n+1/3)+(n+1/3)^2}.
\end{align*}
The signature-three theta functions satisfy many interesting
formulas, including the following cubic relation:
\begin{equation*}
a^3(q)=b^3(q)+c^3(q).
\end{equation*}
Various other properties of $a(q)$, $b(q)$, and $c(q)$ have been
catalogued in \cite{Be5}.

\begin{theorem}\label{hypergeometric theorem} We can reduce the two-dimensional lattice sums to
integrals of hypergeometric functions.

Suppose that $q=e^{-\pi x/\sqrt{12}}$, then
\begin{equation}\label{F(1,1) hypergeometric form}
\begin{split}
\frac{648x}{\pi^2}\frac{F(1,1,1,x^2)}{\left(3+x^2\right)^2}
        &=\left\{\begin{array}{ll}
                     3\tilde{n}\left(3\frac{a( i q)}{b( i q)}\right)
                        +\frac{4}{\sqrt{3}}n_2\left(\frac{b^3( i q)}{a^3( i
                        q)}\right) &\text{if $x\in\left(0,\frac{1}{\sqrt{5}}\right)$},\\
                     3\tilde{n}\left(3\frac{a( i q)}{b( i q)}\right)
                        +\frac{1}{\sqrt{3}}n_2\left(\frac{b^3( i q)}{a^3( i
                        q)}\right) &\text{if $x\in\left(\frac{1}{\sqrt{5}},\sqrt{5}\right)$},\\
                        \frac{1}{\sqrt{3}}n_2\left(\frac{b^3( i q)}{a^3( i
                        q)}\right)  &\text{if $x\in\left(\sqrt{5},\infty\right)$},
                \end{array}\right.
\end{split}
\end{equation}
where
\begin{equation*}
\tilde{n}(k)=\Re\left(\log(k)-\frac{2}{k^3}{_4F_3}\left(\substack{1,1,\frac{4}{3},\frac{5}{3}\\
2,2,2};\frac{27}{k^3}\right)\right),
\end{equation*}
and
\begin{equation*}
n_2(k)=\Im\left(\int_{k}^{1}\frac{{_2F_1}\left(\substack{\frac{1}{3},\frac{2}{3}\\
1}; 1-u\right)}{u}\d u\right).
\end{equation*}
Notice that $\tilde{n}(k)=n(k)$ whenever $|k|$ is sufficiently
large.

Suppose that $q=e^{-\pi x}$ and $x>0$, then
\begin{equation}\label{F(1,2) hypergeometric form}
\frac{16x}{\pi^2}F_{(1,2)}(x)=m\left(\frac{ i
f^4(-q)}{\sqrt{q}f^4\left(-q^4\right)}\right),
\end{equation}
where $m(k)$ is defined in \eqref{m(k) def}.

Suppose that $q=e^{-\pi x/3}$ and $\omega=e^{2\pi  i/3}$, then
\begin{equation}\label{F(1,4)(x) hypergeometric integral}
\begin{split}
\frac{144x}{25\pi^2}F_{(1,4)}(x)
        &=\left\{\begin{array}{ll}
                     m\left(4\frac{\varphi^2(\omega q)}{\varphi^2(-\omega q)}\right)
                        -\frac{3}{4}m_2\left(\frac{\varphi^4(-\omega
                        q)}{\varphi^4(\omega
                        q)}\right) &\text{if $x\in\left(0,\frac{1}{\sqrt{2}}\right)$},\\
                     m\left(4\frac{\varphi^2(\omega q)}{\varphi^2(-\omega q)}\right)
                        +\frac{1}{4}m_2\left(\frac{\varphi^4(-\omega q)}{\varphi^4(\omega
                        q)}\right) &\text{if $x\in\left(\frac{1}{\sqrt{2}},\sqrt{2}\right)$},\\
                        \frac{1}{4} m_2\left(\frac{\varphi^4(-\omega q)}{\varphi^4(\omega
                        q)}\right)  &\text{if $x\in\left(\sqrt{2},\infty\right)$},
                \end{array}\right.
\end{split}
\end{equation}
where $m(k)$ is defined in \eqref{m(k) def}, and
\begin{equation*}
m_2(k):=\Im\left(\int_{k}^{1}
                        \frac{_2F_1\left(\substack{\frac{1}{2},\frac{1}{2}\\ 1}; 1-u\right)}{u}\d
                        u\right).
\end{equation*}

If $q=e^{-\pi x/\sqrt{8}}$, then
\begin{equation}\label{F(2,2)(x) hypergeometric integral}
F_{(2,2)}(x)=\frac{9\pi^2}{256x}\int_{\frac{\varphi^2\left(-q^2\right)}{\varphi^2\left(q^2\right)}}^{1}\frac{\left(3u-1\right)}
{\sqrt{u(1-u)}}{_2F_1}\left(\substack{\frac{1}{2},\frac{1}{2}\\
1}; 1-u^2\right)\d u.
\end{equation}
\end{theorem}
\begin{proof}  The proof of this theorem follows from our ability to
invert theta functions.  First recall the classical inversion
formula for the theta function:
\begin{equation}\label{theta inversion formula}
\varphi^2(q)={_2F_1}\left(\substack{\frac{1}{2},\frac{1}{2}\\ 1};
1-\frac{\varphi^4(-q)}{\varphi^4(q)}\right),
\end{equation}
which holds whenever $q\in(-1,1)$ \cite{Be3}.  If we use the
notation $\alpha=1-\varphi^4(-q)/\varphi^4(q)$ and $z=\varphi^2(q)$,
then many different theta functions can be expressed in terms of
these two parameters.  The following identities are true whenever
$|q|<1$:
\begin{align*}
\varphi(q)=&\sqrt{z},\\
\varphi(-q)=&(1-\alpha)^{1/4}\sqrt{z},\\
\varphi\left(q^2\right)=&\left(1+\sqrt{1-\alpha}\right)^{1/2}\sqrt{\frac{z}{2}},\\
\psi(-q)=&q^{-1/8}\left\{\alpha(1-\alpha)\right\}^{1/8}\sqrt{\frac{z}{2}},\\
\psi(q)=&q^{-1/8}\alpha^{1/8}\sqrt{\frac{z}{2}},
\end{align*}
and it is also well known that
\begin{equation*}
\frac{\d \alpha}{\d q}=\frac{\alpha(1-\alpha)z^2}{q}.
\end{equation*}
Both \eqref{F(1,2) hypergeometric form} and \eqref{F(2,2)(x)
hypergeometric integral} follow from applying these
parameterizations to equations \eqref{F(1,2)(x) q-series integral}
and \eqref{F(2,2)(x) q-series integral} respectively.

Since equation \eqref{theta inversion formula} does not hold in the
entire open unit disk, we will need to generalize that result. First
notice that $z$ satisfies the hypergeometric differential equation
with respect to $\alpha$:
\begin{equation}\label{hypergeometric diff equation}
\alpha(1-\alpha)\frac{\d^2 z}{\d\alpha^2}+(1-2\alpha)\frac{\d
z}{\d\alpha}-\frac{z}{4}=0.
\end{equation}
We can use the relation
$\frac{\d}{\d\alpha}=\frac{1}{\frac{\d\alpha}{\d
q}}\times\frac{\d}{\d q}$, to show that \eqref{hypergeometric diff
equation} holds (excluding possible poles) for $|q|<1$.  The most
general solution of this differential equation has the form
\begin{equation*}
z=C {_2F_1}\left(\substack{\frac{1}{2},\frac{1}{2}\\ 1};
\alpha\right)+D {_2F_1}\left(\substack{\frac{1}{2},\frac{1}{2}\\
1} ; 1-\alpha\right),
\end{equation*}
where $C$ and $D$ are undetermined constants.  When $q$ lies in a
neighborhood of zero, \eqref{theta inversion formula} shows that
$(C,D)=(1,0)$.  We can analytically continue that solution to a
larger connected $q$-domain, provided that $\alpha$ (and $1-\alpha$
if $D\not=0$) does not intersect the line
$[1,\infty)$.  In particular, the function ${_2F_1}\left(\substack{\frac{1}{2},\frac{1}{2}\\
1}; \alpha\right)$ has a branch cut running along $[1,\infty)$.

If we consider values of $q\in(0,\omega)$ with $\omega=e^{2\pi
i/3}$, then $\alpha$ crosses $[1,\infty)$ at the point $q=\omega
e^{-\pi\sqrt{2}/3}$. Similarly, $1-\alpha$ intersects the branch cut
at $q=\omega e^{-\pi/3\sqrt{2}}$. It follows that we will have to
solve the hypergeometric differential equation separately on each of
the three line segments. If $u=\omega e^{-\pi x/3}$, then
\begin{equation}\label{z inversion on a ray}
\begin{split}
\varphi^2(u)    &=\left\{\begin{array}{ll}
                      -3{_2F_1}\left(\substack{\frac{1}{2},\frac{1}{2}\\ 1};1-\frac{\varphi^4(-u)}{\varphi^4(u)}\right)+2 i{_2F_1}\left(\substack{\frac{1}{2},\frac{1}{2}\\ 1};\frac{\varphi^4(-u)}{\varphi^4(u)}\right) &\text{if $x\in\left(0,\frac{1}{\sqrt{2}}\right)$},\\
                      {_2F_1}\left(\substack{\frac{1}{2},\frac{1}{2}\\
                        1};1-\frac{\varphi^4(-u)}{\varphi^4(u)}\right)+2 i
                        {_2F_1}\left(\substack{\frac{1}{2},\frac{1}{2}\\ 1};\frac{\varphi^4(-u)}{\varphi^4(u)}\right) &\text{if $x\in\left(\frac{1}{\sqrt{2}},\sqrt{2}\right)$},\\
                        {_2F_1}\left(\substack{\frac{1}{2},\frac{1}{2}\\ 1};1-\frac{\varphi^4(-u)}{\varphi^4(u)}\right)  &\text{if
                        $x\in\left(\sqrt{2},\infty\right)$}.
                \end{array}\right.
\end{split}
\end{equation}
The coefficients in \eqref{z inversion on a ray} can be verified
from the fact that $\varphi^2(u)$ is analytic when $x\in
(0,\infty)$. For example, we can check the continuity of the
right-hand side of \eqref{z inversion on a ray} by letting
$u\rightarrow\omega e^{-\pi\sqrt{2}/3}$. In that case
$\alpha=1-\frac{\varphi^4(-u)}{\varphi^4(u)}\approx 5.828\dots$, and
we have:
\begin{equation*}
\begin{split}
0=&\varphi^2\left(\omega e^{-\pi\frac{\sqrt{2}+0}{3}}\right)-\varphi^2\left(\omega e^{-\pi\frac{\sqrt{2}-0}{3}}\right)\\
=&{_2F_1}\left(\substack{\frac{1}{2},\frac{1}{2}\\ 1};\alpha+ i 0\right)-{_2F_1}\left(\substack{\frac{1}{2},\frac{1}{2}\\ 1};\alpha- i 0\right)-2 i{_2F_1}\left(\substack{\frac{1}{2},\frac{1}{2}\\
1};1-\alpha\right).
\end{split}
\end{equation*}
This vanishing of this last expression follows from basic properties
of the hypergeometric function (see problem $1$ on page 276 of
\cite{Le}), and therefore the right-hand side of \eqref{z inversion
on a ray} is indeed continuous at $x=\sqrt{2}$. In practice, we
simply discovered \eqref{z inversion on a ray} numerically.

We will use the theory of signature-three theta functions to prove
equation \eqref{F(1,1) hypergeometric form}. Recall that $c(q)$ can
be expressed as an infinite product:
\begin{equation*}
\frac{c^3(q)}{27 q}=\frac{f^9(-q^3)}{f^3(-q)},
\end{equation*}
and that the signature-three theta functions obey a differentiation
formula:
\begin{equation*}
\frac{c^3(q)}{q}=\frac{a(q)}{1-\frac{c^3(q)}{a^3(q)}}\frac{\d}{\d
q}\left(\frac{c^3(q)}{a^3(q)}\right).
\end{equation*}
It follows immediately that equation \eqref{F(1,1)(x) q-series
integral} reduces to
\begin{equation}\label{F(1,1)(x) q integral reduced}
\frac{F(1,1,1,x^2)}{(3+x^2)^2}=\frac{\pi^2}{648\sqrt{3}x}\Im\left(\int_{0}^{
i q}\frac{a(u)}{1-\frac{c^3(u)}{a^3(u)}}\frac{\d}{\d
u}\left(\frac{c^3(u)}{a^3(u)}\right)\d u\right).
\end{equation}
Next recall that for $|u|$ sufficiently small:
\begin{equation}\label{a(q) inversion formula}
a(u)={_2F_1}\left(\substack{\frac{1}{3}, \frac{2}{3}\\
1};\frac{c^3(u)}{a^3(u)}\right).
\end{equation}
In order to apply \eqref{a(q) inversion formula} to our integral, we
will need to establish a generalized inversion formula which holds
for $u\in (0, i)$.  The reasoning closely follows the proof of
\eqref{z inversion on a ray}, except that
$c^3(u)/a^3(u)\in[1,\infty)$ when $u= i e^{-\pi\sqrt{5/12}}$, and
$1-c^3(u)/a^3(u)\in[1,\infty)$ when $u= i e^{-\pi/\sqrt{60}}$.
Suppose that $u= i e^{-\pi x/\sqrt{12}}$, then we obtain
\begin{equation}\label{a(q) inversion on a ray}
\begin{split}
a(u)        &=\left\{\begin{array}{ll}
                      4{_2F_1}\left(\substack{\frac{1}{3},\frac{2}{3}\\ 1};\frac{c^3(u)}{a^3(u)}\right)
                      +\sqrt{3} i{_2F_1}\left(\substack{\frac{1}{3},\frac{2}{3}\\ 1};1-\frac{c^3(u)}{a^3(u)}\right) &\text{if $x\in\left(0,\frac{1}{\sqrt{5}}\right)$},\\
                      2{_2F_1}\left(\substack{\frac{1}{3},\frac{2}{3}\\ 1};\frac{c^3(u)}{a^3(u)}\right)
                      +\sqrt{3} i{_2F_1}\left(\substack{\frac{1}{3},\frac{2}{3}\\ 1};1-\frac{c^3(u)}{a^3(u)}\right) &\text{if $x\in\left(\frac{1}{\sqrt{5}},\sqrt{5}\right)$},\\
                      {_2F_1}\left(\substack{\frac{1}{3},\frac{2}{3}\\ 1};\frac{c^3(u)}{a^3(u)}\right) &\text{if $x\in\left(\sqrt{5},\infty\right)$}.
                \end{array}\right.
\end{split}
\end{equation}
Finally, \eqref{F(1,1) hypergeometric form} follows from
substituting \eqref{a(q) inversion on a ray} into \eqref{F(1,1)(x) q
integral reduced} and simplifying. $\blacksquare$
\end{proof}

Finally, we will conclude this section by summarizing the formulas
that follow from setting $x=1$ in Theorem \ref{hypergeometric
theorem}.

\begin{corollary}\label{final reduction cor} The following identities are true:
\begin{align}
\frac{27}{2\pi^2}F(1,1)=&n\left(3\sqrt[3]{2}\right),\\
\frac{16}{\pi^2}F(1,2)=&m\left(4 i\right),\\
\frac{144}{25\pi^2}F(1,4)=&m\left(\frac{4}{\theta}\right)
+\frac{1}{4}\Im\left(\int_{\theta^2}^{1}\frac{{_2F_1}\left(\substack{\frac{1}{2},\frac{1}{2}\\
1}; 1-u\right)}{u}\d u\right)\label{F(1,4) final formula},\\
\frac{256}{9\pi^2}F(2,2)=&\int_{\sqrt{2}-1}^{1}\frac{3u-1}{\sqrt{u(1-u)}}{_2F_1}\left(\substack{\frac{1}{2},\frac{1}{2}\\
1}; 1-u^2\right)\d u,\label{F(2,2) final formula}
\end{align}
where $\theta=\frac{(4-2t-2t^2+t^3)}{4\sqrt{2}}$, and $t=-
i\sqrt[4]{12}$.
\end{corollary}
Notice that equation \eqref{F(1,4) final formula} involves Meijer's
$G$-function disguised as a hypergeometric integral:
\begin{equation*}
\Im\left(\int_{k}^{1}\frac{{_2F_1}\left(\substack{\frac{1}{2},\frac{1}{2}\\
1};1-u\right)}{u}\d
u\right)=\frac{1}{\pi^2}\Im\left({G_{3,3}^{3,2}}\left(k\big\vert\substack{\frac{1}{2},\frac{1}{2},1\\
0,0,0}\right)\right).
\end{equation*}
This identity probably rules out the possibility of expressing
$F(1,4)$ as a Mahler measure, and it also indicates that any
explicit formula for $F(b,c)$ should reduce to Meijer $G$-functions
in certain instances.

\section{Formulas for $F(1,1,2,4)$, $F(1,2,4,4)$, and additional explicit examples}\label{section: F(1,1,2,4) and F(1,2,4,4)}\init

In this section, we will present several additional formulas for
$F(a,b,c,d)$.  Since the ideas are the same as in the previous
section, we will only sketch brief details of each proof.

\begin{theorem} The following identities are true:
\begin{align}
\frac{144\sqrt{2}}{121\pi^2}F(1,2,4,4)&=\frac{1}{3}n\left(3\left(\frac{4}{2+17\sqrt{2}-9\sqrt{6}}\right)^{1/3}\right)-\frac{1}{6}n\left(3\left(\frac{4}{2+\sqrt{2}}\right)^{1/3}\right),\label{F(1,1,2,4) identity}\\
\frac{27}{\sqrt{2}\pi^2}F(1,1,2,4)&=n\left(3\left(\frac{4}{2-17\sqrt{2}+9\sqrt{6}}\right)^{1/3}\right)-n\left(3\left(\frac{4}{2-\sqrt{2}}\right)^{1/3}\right).\label{F(1,2,4,4)
identity}
\end{align}
\end{theorem}
\begin{proof} The proof of this theorem requires the following identities:
\begin{align*}
e_1^2 e_2 e_4=&\left(\frac{e_1^2
e_4^2}{e_2}\right)\left(\frac{e_2^2}{e_4}\right)=\sum_{n,k=-\infty}^{\infty}(-1)^k(3n+1)q^{\frac{(3n+1)^2+6k^2}{3}},\\
e_1 e_2 e_4^2=&\left(\frac{e_1^2
e_4^2}{e_2}\right)\left(\frac{e_2^2}{e_1}\right)=\sum_{\substack{n=-\infty\\
k=0}}^{\infty}(3n+1)q^{\frac{8(3n+1)^2+3(2k+1)^2}{24}}.
\end{align*}
By an argument similar to the one in Proposition \ref{q-series
proposition}, we can prove that
\begin{equation*}
\begin{split}
F(1,1,2,4)=&\sum_{n,k=-\infty}^{\infty}\frac{(-1)^k(3n+1)}{\left((3n+1)^2+6k^2\right)^2}\\
=&-\frac{\sqrt{2}\pi^2}{81}\left(\log(q^4)+9\sum_{n=1}^{\infty}n\chi_{-3}(n)\log\left(1+q^{4n}\right)\right),
\end{split}
\end{equation*}
and
\begin{equation*}
\begin{split}
F(1,2,4,4)=&\frac{121}{2}\sum_{n,k=-\infty}^{\infty}\frac{(3n+1)}{\left(8(3n+1)^2+3(2k+1)^2\right)^2}\\
=&\frac{121\pi^2}{288\sqrt{2}}\sum_{n=1}^{\infty}n\chi_{-3}(n)\log\left(\frac{1+q^n}{1-q^n}\right),
\end{split}
\end{equation*}
 where $q=e^{-\pi/\sqrt{6}}$.  Both of these $q$-series can be evaluated in terms of $n(k)$ by the results
 of Rodriguez-Villegas (see \cite{LR} for details). $\blacksquare$
\end{proof}

If we consider additional examples involving $F(1,1,1,x^2)$, then we
can establish many interesting formulas by setting $x=\sqrt{a/b}$.
For instance, if $x\in\left\{1,\frac{1}{\sqrt{5}}\right\}$, we have
$\Im\left(b^3( i q)/a^3( i q)\right)=0$, and as a result it is
possible to show that the $n_2$ term in \eqref{F(1,1) hypergeometric
form} vanishes.

\begin{theorem}Let $\phi=\frac{1+\sqrt{5}}{2}$, $y=\sqrt[3]{2}e^{\pi i/3}$, and suppose that $z\approx-.58+.56 i$ is a root of the equation $(z^2+3z+1)^3-2(z^6+1)=0$, then
\begin{align}
F(1,1)=&\frac{2\pi^2}{27}n\left(3\sqrt[3]{2}\right),\label{F(1,1)(1) to ntilde}\\
F\left(1,1,1,\frac{1}{5}\right)=&\frac{32\pi^2}{135\sqrt{5}}\tilde{n}\left(\frac{3}{\sqrt[3]{\phi}}\right),\label{F(1,1)sqrt5}\\
F\left(1,1,1,5\right)\stackrel{?}{=}&\frac{4\pi^2}{27\sqrt{5}}n\left(-3\sqrt[3]{\phi}\right),\label{F(1,1)sqrt5b}\\
F(1,1,1,9)=&
\frac{2\pi^2}{27\sqrt{3}}n_2\left(\frac{9}{2(1+y+3y^2)}\right),\\
F(1,1,1,25)=&
\frac{98\pi^2}{405\sqrt{3}}n_2\left(\frac{1}{1+z^6}\right).
\end{align}
\end{theorem}

Notice that the result in the introduction, \eqref{intro lattic
example}, follows from combining \eqref{F(1,1)sqrt5} with
\eqref{n(k) 4F3 to 3F2}.  We have also stated \eqref{F(1,1)sqrt5b}
as a conjecture, because the result presumably follows from
\eqref{F(1,1) hypergeometric form} when $x\rightarrow{\sqrt{5}}$,
but we still have not found a completely rigorous proof of that
fact.
 Finally, we will conclude this section with two additional identities.
If we set $q=e^{-\frac{\pi}{2\sqrt{3}}}$ and $\omega=e^{\pi i/4}$,
then
\begin{align}
\frac{24^2\sqrt{2}}{19^2
\pi^2}F(1,2,8,8)=&\sum_{n=1}^{\infty}n\chi_{-3}(n)\log\left|\frac{1+\omega
q^n}{1-\omega q^n}\right|,\label{F(1,2,8,8) formula}\\
\frac{72\sqrt{2}}{49\pi^2}\F(1,1,4,8)=&\sum_{n=1}^{\infty}n\chi_{-3}(n)\log\left|\frac{1+\omega
q^{2n}}{1-\omega q^{2n}}\right|.\label{F(1,1,2,8) formula}
\end{align}
Unfortunately, we have not been able to reduce these last two
$q$-series to known functions.  These $q$-expansions follow from
combining $e_1 e_2 e_8^2=\left(\frac{e_2^2
e_8^2}{e_4}\right)\left(\frac{e_1 e_4}{e_2}\right)$ and $e_1^2 e_4
e_8=\left(\frac{e_1^2 e_4^2}{e_2}\right)\left(\frac{e_2
e_8}{e_4}\right)$, with equations \eqref{weight 1/2 series} and
\eqref{weight 3/2 series}.

\section{Remarks on $F(1,3)$, $F(2,9)$, $F(5,9)$ and higher values of
$F(a,b,c,d)$}\label{section: F(1,3) section}\init

In this section, we will briefly demonstrate how to reduce several
higher values of $F(a,b,c,d)$ to known functions.  Our proof of the
$F(1,3)$ formula will be instructive. Recall that Rodriguez-Villegas
showed that
\begin{equation}\label{F(1,3) RV formula}
\frac{4\pi^2}{81}n(-6) 
=\Re\left(\frac{1}{2}\sum_{\substack{m,n\in\mathbb{Z}\\
(m,n)\not=(0,0)}}\frac{\chi_{-3}(n)} {\left(3\left(\frac{1+
i\sqrt{3}}{2}\right)m+n\right)^2\left(3\left(\frac{1-
i\sqrt{3}}{2}\right)m+n\right)}\right),
\end{equation}
and then used Deuring's theorem to equate this Eisenstein series to
the $L$ series of a CM elliptic curve of conductor $27$
\cite[p.~32]{RV}. A different proof could have been constructed from
numerically observing that
\begin{equation}\label{F(1,3) essential q-series identity}
\begin{split}
q\prod_{n=1}^{\infty}&\left(1-q^{3n}\right)^2\left(1-q^{9n}\right)^2\\
&=\frac{1}{4}\sum_{j=1}^{2}\chi_{-3}(j)
\sum_{n,m=-\infty}^{\infty}\left((6m+j)+3(6n+j)\right)q^{\frac{(6m+j)^2+3(6n+j)^2}{4}}.
\end{split}
\end{equation}
The modularity theorem implies that $e_3^2 e_9^2$ is associated to
the correct elliptic curve, hence the Mellin transform of the
left-hand side of equation \eqref{F(1,3) essential q-series
identity} will equal $L(E,s)$.  Since the Mellin transform (at
$s=2$) of the right-hand side trivially equals the right-hand side
of equation \eqref{F(1,3) RV formula}, it just remains to prove
\eqref{F(1,3) essential q-series identity}.  By applying limiting
cases of the triple and quintuple product identities, we can show
that equation \eqref{F(1,3) essential q-series identity} is
equivalent to an identity between eta functions:
\begin{equation}\label{F(1,3) equivalent eta function identity}
\begin{split}
4e_3^2e_9^2
=&\left(\frac{e_6^5e_{36}e_{54}^2}{e_{12}^2e_{18}e_{108}}\right)
+3\left(\frac{e_{12} e_{18}^7}{e_{6}e_{36}^3}\right)
-2\left(\frac{e_3^2e_{12}^2e_{18}^2e_{27}e_{108}}
{e_6e_9e_{36}e_{54}}\right)
-6\left(\frac{e_6^2e_9^3e_{36}^3}{e_3e_{12}e_{18}^2}\right).
\end{split}
\end{equation}
Identities between modular forms, such as \eqref{F(1,3) equivalent
eta function identity}, are usually established by checking that
both sides of a formula have the same McLauren series expansion for
sufficiently many terms.  If we use the inversion formula for the
eta function, then \eqref{F(1,3) equivalent eta function identity}
can also be viewed as an example of a mixed modular equation
\cite{Be3}. Of course, the main difficulty of extending this type of
approach, is to actually find the necessary modular equations. Since
the first version of this paper appeared, it has become clear that
this method also applies to $F(2,9)$, $F(5,9)$ and $F(4,7,7,28)$
\cite{RY}.


\section{New Mahler measures, and the $L$-series of an irrational modular form}\label{section: more mahlers}\init

In this section we will use values of class invariants to deduce
some explicit formulas for Mahler measures.  Recall that if
$q=e^{-\pi\sqrt{m}}$, then the class invariants are defined by
\begin{align*}
g_m:=&2^{-1/4}q^{-1/24}\left(q;q^2\right)_{\infty},&
G_m:=&2^{-1/4}q^{-1/24}\left(-q;q^2\right)_{\infty}.
\end{align*}
It is a classical fact that $G_m$ and $g_m$ are algebraic numbers
whenever $m\in\mathbb{Q}$, and that they satisfy the following
algebraic relation:
\begin{equation}\label{Gm and gm algebraic relation}
\left(g_m G_m\right)^8\left(G_m^8-g_m^8\right)=\frac{1}{4}.
\end{equation}
Since most tables only contain values of $g_m$ when $m$ is even, and
$G_m$ when $m$ is odd, our calculations will require \eqref{Gm and
gm algebraic relation}.

\begin{theorem}\label{Theorem: addition F(1,2)(m) cases} Suppose that $m\in\mathbb{N}$, then
\begin{equation}
m\left(8 i g_m^8
G_m^4\right)=\frac{16\sqrt{m}}{\pi^2}\sum_{n=1}^{\infty}\frac{b_n}{n^2},
\end{equation}
where
\begin{equation*}
\sum_{n=1}^{\infty}b_n q^n=\frac{e_8^3 e_{4m}^2}{e_{8m}}.
\end{equation*}
The following table gives evaluations of $8g_m^8 G_m^4$, and states
whether or not $b_n$ is multiplicative:
\begin{equation*}
    \begin{tabular}{|c|c|c|p{6 in}|}
        \hline
        $m$ &  $8g_m^8 G_m^4$ &   $b_n$ multiplicative?\\
        \hline
        $1$ & $4$ &   Yes\\
        $2$ & $4\sqrt{2+2\sqrt{2}}$ &   No\\
        $3$ & $4\left(2+\sqrt{3}\right)$ &   No\\
        $7$ & $4\left(8+3\sqrt{7}\right)$ & No\\
        $9$ & $4\left(7+4\sqrt[4]{12}+2\sqrt[4]{12^2}+\sqrt[4]{12^{3}}\right)$ &  No\\
        $15$ & $4\left(28+16\sqrt{3}+12\sqrt{5}+7\sqrt{15}\right)$ &  No\\
        \hline
    \end{tabular}
\end{equation*}
\end{theorem}
\begin{proof}
If $m\in\mathbb{N}$, then we can use the definition of
$F_{(1,2)}(x)$ to show that
\begin{equation*}
F_{(1,2)}(\sqrt{m})=\sum_{n=1}^{\infty}\frac{b_n}{n^2},
\end{equation*}
where $b_n$ has the stated generating function. Furthermore,
equation \eqref{F(1,2)(x) as a Mahler measure} reduces to
\begin{equation}
F_{(1,2)}\left(\sqrt{m}\right)=\frac{\pi^2}{16\sqrt{m}}m\left(8 i
g_m^8 G_m^4\right).
\end{equation}
Therefore, we can obtain Mahler measure formulas by appealing to
tables of class invariants \cite{Be5}. In order to check the values
of $8g_m^8 G_m^4$, we can solve \eqref{Gm and gm algebraic relation}
to show that
\begin{align*}
8g_m^8 G_m^4&=4\left(G_m^{12}+\sqrt{G_m^{24}-1}\right)\\
 &=4\sqrt{2}
g_m^{6}\sqrt{g_m^{12}+\sqrt{g_m^{24}+1}}.
\end{align*}
For example, since $G_1=1$, it follows that $8g_1^8 G_1^4=4$.  While
this type of argument frequently leads to nested radicals, many of
those identities simplify with sufficient effort. $\blacksquare$
\end{proof}

A cursory inspection of Theorem \ref{Theorem: addition F(1,2)(m)
cases}, reveals only one instance where $b_n$ is multiplicative. In
particular, when $m=1$ the generating function for $b_n$ reduces to
$e_4^2 e_8^2$.  This cusp form is associated to a conductor $36$
elliptic curve with complex multiplication \cite{Ono}.  Thus, we
have given a new proof of the formula for $m(4i)$ \cite{RV}.  When
$m=2$, a brief computation reveals that $b_3=b_{11}=0$, but
$b_{33}=-8$, and therefore $e_8^5/e_{16}$ is not multiplicative.
Fortunately, Somos has pointed out that $e_8^5/e_{16}$ is the real
part of a multiplicative cusp form in $\mathbb{Q}(i\sqrt{2})$:
\begin{equation*}
\begin{split}
g(q)&=\frac{e_8^5}{e_{16}}+2i\sqrt{2}\frac{e_{16}^5}{e_{8}}\\
&=q+2i\sqrt{2} q^3-5q^9+\dots
\end{split}
\end{equation*}
It is not difficult to find a formula for $L(e_{16}^{5}e_8^{-1},2)$,
and therefore we can show that
\begin{equation}\label{L(g,2)}
\frac{16\sqrt{2}}{\pi^2}L(g,2)=m\left(4i\sqrt{2\sqrt{2}+2}\right)+\frac{i}{\sqrt{2}}m\left(4\sqrt{2\sqrt{2}-2}\right).
\end{equation}
It would probably be interesting to determine if $g(q)$ holds some
special arithmetic significance.
%

\section{Higher polylogarithms and conclusion}
\label{section: conclusion} \init

We will conclude the paper, by showing that the method from Section
\ref{section:hypergeometric reductions} can be used to produce
identities for elliptic polylogarithms. Let us briefly reexamine
Proposition \ref{q-series proposition}.  If we had used involutions
for weight $1/2$ theta functions, rather than involutions for weight
$3/2$ theta functions, we would have obtained formulas including
\begin{equation*}
\frac{F(1,1,1,x)}{(3+x)^2}=\frac{\pi}{36\sqrt{x}}\sum_{n=1}^{\infty}\chi(n)
D\left(i e^{-\frac{\pi n}{\sqrt{12x}}}\right),
\end{equation*}
where $D(z)=\Im\left(\Li_2(z)+\log|z|\log(1-z)\right)$ is the
Bloch-Wigner dilogarithm, and $\chi(n)$ is a character modulo $12$,
with $\chi(1)=\chi(11)=1$, and $\chi(5)=\chi(7)=-1$.  It is also
interesting to note that the following conjecture of
Rodriguez-Villegas \cite{Fi}:
\begin{equation*}
\m\left(1+x_1+x_2+x_3+x_4\right)\stackrel{?}{=}\frac{675\sqrt{15}}{16\pi^5}L(f,4),
\end{equation*}
where $f(q)=e_3^3e_5^3+e_1^3e_{15}^3$, can be reformulated using
such an argument.  By the method of Proposition \ref{q-series
proposition}, we have
 \begin{equation*}
 L(f,4)=-\frac{128\pi}{15^3}\sum_{n=1}^{\infty}\chi_{-4}(n)\left(9R\left(i
 q^{n}\right)+R\left(i
 q^{3n}\right)\right),
 \end{equation*}
 where
 $R(z)=\Im\left(\log|z|\Li_4(z)-\log^2|z|\Li_3(z)+\frac{\log^3|z|}{3}\Li_2(z)\right)$,
 and $q=e^{-\pi\sqrt{15}/6}$.

\bigskip


\begin{acknowledgements}

The author would like to thank David Boyd, Wadim Zudilin, Larry
Glasser, Bruce Berndt, Michael Somos, and Anton Mellit for their
useful comments and encouragement.  Additionally, the author thanks
the Max Planck Institute of Mathematics for their hospitality.
\end{acknowledgements}

\end{document}